\documentclass{amsart}
\usepackage{amsmath,amssymb,amscd}

\newcommand{\entirea}{\varprojlim H_{\varepsilon}^{\text{{\rm ev}}}(\mathfrak{A}_n)}
\newcommand{\entireb}{H_{\varepsilon}^{\text{{\rm ev}}}(\mathfrak{A})}
\newcommand{\retsu}{a_0, \cdots , a_{2k}}
\newcommand{\fre}{\mathfrak{A}_n}
\newcommand{\tori}{(T^2_{\theta})^\infty}
\newcommand{\C}{\mathbb{C}}
\newcommand{\Z}{\mathbb{Z}}
\newcommand{\R}{\mathbb{R}}
\newcommand{\nipai}{2\pi i}
\newcommand{\ide}{\mathfrak{I}}
\newcommand{\e}{\varepsilon}
\newcommand{\ct}{C^\infty(T)}
\newcommand{\ore}{\mathfrak{A}}

\DeclareMathOperator{\supp}{supp}
\DeclareMathOperator{\ran}{Ran}

\newtheorem{lemma}{Lemma}[section]
\newtheorem{proposition}{Proposition}[section]
\newtheorem{definition}{Definition}[section]
\newtheorem{theorem}{Theorem}[section]

\begin{document}

\title{The Entire Cyclic Cohomology of Noncommutative 2-tori}
\author{Katsutoshi Kawashima$^{1}$}\footnote{
           Department of Mathematics and Information Sciences, 
           Tokyo Metropolitan University \\
           e-mail:naito-katsutoshi@ed.tmu.ac.jp  \quad The former family name was NAITO.}

\maketitle

%\begin{center}
%{\it Department of Mathematics and Information Sciences \\
%Tokyo Metropolitan University}
%\end{center}

\begin{abstract}
Our aim in this paper is to compute the entire cyclic cohomology of noncommutative 2-tori. First of all, we clarify their algebraic structure of noncommutative 2-tori as a $F^*$-algebra, according to the idea of Elliott-Evans. Actually, they are the $F^*$-inductive limit of subhomogeneous $F^*$-algebras. Using such a result, we compute their entire cyclic cohomology, which is isomorphic to their periodic one as a complex vector space.
\end{abstract}

\section{Introduction}

Elliott and Evans \cite{ee} show that the irrational rotation $C^*$-algebras (or noncommutative 2-tori) $T^2_\theta$ are isomorphic to certain inductive limits, which are now called AT-algebras, 
\[ \varinjlim (C(T)\otimes (M_{q_{2n}}(\C)\oplus M_{q_{2n-1}}(\C)), \pi_n). \]

To compute the entire cyclic cohomology of their smooth parts $\tori$, we need to know their algebraic structure. In this paper, we elaborate Elliott and Evans' result cited above, and show that $\tori$ are isomorphic to inductive limits 
\[ \varinjlim (\ct\otimes (M_{q_{2n}}(\C)\oplus M_{q_{2n-1}}(\C)), \pi_n^\infty) \]
as Fr\'echet $^*$-algebras (or $F^*$-algebras). Using this fact, we can compute their entire cyclic cohomology quite easily.

In Sect.\ref{pre}, we prepare the notations needed for $\tori$ and review the definition of entire cyclic cohomology. In Sect.\ref{str}, we determine the algebraic structure of $\tori$ by using appropriate smooth functions to construct projections based on Connes \cite{co2} instead of the original ones due to Rieffel \cite{r}. In Sect.\ref{fre}, it is shown that the functor of entire cyclic cohomology $H_\e^*$ is continuous in some sense. More precisely,
\[ H_\e^* (\varinjlim \mathfrak{A}_n )\simeq 
\varprojlim H_\e^* (\mathfrak{A}_n) \]
(cf. Meyer \cite{meyer}), where the right hand side means the projective limit of $H_\e^*(\fre)$ which will be defined in the same section.

Our main result is stated in Sect.\ref{main}.

\section{Preliminaries}\label{pre}

First of all, we define some notations for our discussion in this section. 

Given an irrational number $\theta$, let us treat the noncommutative 2-tori $\tori$ generated by two unitaries $u, v$  with reltaion
\[ uv=e^{\nipai\theta}vu \]
as a Fr\'echet *-algebra (or $F^*$-algebra).
In some cases, we regard each element of $\tori$ as an operator on the Hilbert space $L^2(T)$ of the square integrable complex valued functions on the 1-torus $T$. For instance, 
\[ (uf)(t)=tf(t), \quad (vf)(t)=f(e^{-\nipai\theta}t) \]
for $f\in L^2(T), \, t\in T$. 

There is a smooth action $\alpha$ of $T^2$ on $\tori$ defined by 
\[ \alpha_{t, s}(u)=tu, \quad \alpha_{t, s}(v)=sv \]
for $t, s\in T$. Moreover, we have the two *-derivations $\delta_1, \delta_2$ on $\tori$ associated with $\alpha$  satisfying
\[ \delta_1(u)=iu, \quad \delta_2(u)=0, \quad \delta_1(v)=0, \quad \delta_2(v)=iv. \]
Using these derivations, we define seminorms $\|\cdot\|_{k, l}$ on $\tori$ by
\[ \|x\|_{k, l}=\|\delta_1^k\circ\delta_2^l (x)\|, \]
where $\|\cdot\|$ is the usual $C^*$-norm on $T^2_\theta$.

Here, we briefly review the definition of entire cyclic cohomology. For any unital $F^*$-algebra $\mathfrak{A}$ and any integer $n\geq 0$, we put $C^n$ be the set of all $(n+1)-$linear functionals on $\mathfrak{A}$. For $n<0$, let $C^n =\{ 0\}$. Moreover, we define
\begin{align*}
C^{\rm{ev}}&=\{(\varphi_{2n})_n \,|\, \varphi_{2n}\in C^{2n} \,(n\geq 0) \}, \\
C^{\rm{od}}&=\{(\varphi_{2n+1})_n \,|\, \varphi_{2n+1}\in C^{2n+1} \,(n\geq 0) \}.
\end{align*}

We call $(\varphi_{2n})$ an entire even cochain if for each bounded subset $\Sigma \subset \mathfrak{A}$, we can find a constant $C>0$ such that
\[ \left|\varphi_{2n}(a_0, \dots , a_{2n})\right|\leq C\cdot n! \]
for all $n\geq 1$ and $a_j\in \Sigma$.
 In odd case, we define entire odd cochains by the same way as in even case. We denote by $C_\e^{\rm{ev}}$ (resp. $C_\e^{\rm{od}}$) the set of all entire even (resp. odd) cochains. Then we define the entire cyclic cohomology of $\mathfrak{A}$ by the cohomology of the short complex
 \[ C_\e^{\rm{ev}} \overset{\partial}{\underset{\partial}\leftrightarrows} C_\e^{\rm{od}}, \]
 where $\partial$ are certain derivativions defined by Connes \cite{co}.

\section{$\tori$ is a Fr\'echet Inductive Limit}\label{str}

In this section, we prove the key lemma which states that noncommutative 2-tori $\tori$ as $F^*$-algebras are isomorphic to inductive limits
\[ \varinjlim (C^\infty (T)\otimes (M_{q_{2n}}(\C)\oplus M_{q_{2n-1}}(\C)), \pi^\infty_n), \]
where the sequence $\{q_{2n-1}\}_n$ appears in the continued fraction expansion of $\theta$.

Let $\begin{pmatrix}
p' & p \\
q' & q
\end{pmatrix} \in SL(2, \Z)$
with $p/q<\theta<p'/q', \, q>0$ and $q'>0$ for each fixed $\theta \in (0, 1)$. We write $\beta=p'-q'\theta, \beta'=q\theta-p$. First of all, we construct two projections $e_\beta$ and $e_{\beta '}$ in $\tori$ with traces $\beta$ and $\beta'$ respectively using the functions $f_\beta$ and $g_\beta$ defined below. We regard the $1$-torus $T$ as the interval $[0, 1]$. Since 
$\begin{pmatrix}
p' & p \\
q' & q
\end{pmatrix}\in SL(2, \Z)$, we note that $q\beta+q'\beta'=1$. In particular, we have $0<\beta<1/q, 0<\beta'<1/q'$. When $\beta\geq 1/2q$, we put 
\begin{align*}
f_1(x)&=e^{-\alpha/x} &
f_2(x)&=1-f_1(1/q-\beta-x) \\
f_3(x)&=f_2(1/q-x) & f_4(x)&=f_1(1/q-x),
\end{align*}
where $\alpha=(1/q-\beta)\log \sqrt{2}$. Using the functions described above, we define the functions $f, g$ defined by
\begin{align*} f_\beta(x)&=
\begin{cases}
f_1(x) & (0\leq x \leq 1/2q-\beta/2) \\
f_2(x) & (1/2q-\beta/2 \leq x \leq 1/q-\beta) \\
1 & (1/q-\beta \leq x \leq \beta) \\
f_3(x) & (\beta \leq x \leq \beta/2+1/2q) \\
f_4(x) & (\beta/2+1/2q \leq x \leq 1/q) \\
0 & (1/q \leq x <1),
\end{cases} 
\\
g_\beta(x)&=\chi_{[\beta, 1/q]}(x)\sqrt{f(x)-f(x)^2}, 
\end{align*}
where $\chi$ stands for the characteristic function. In the case when $\beta <1/2q$, we put 
\begin{align*}
f_1(x)&=e^{-\alpha'/x} & f_2(x)&=1-f_1(1/q-\beta-x) \\
f_3(x)&=f_2(\beta-x) & f_4(x)&=f_1(\beta-x),
\end{align*}
where $\alpha'=\beta\log \sqrt{2}$, and define
\begin{align*} f_\beta(x)&=
\begin{cases}
f_1(x) & (1/2q-\beta \leq x \leq 1/2q-\beta/2) \\
f_2(x) & (1/2q-\beta/2 \leq x \leq 1/2q) \\
f_3(x)& (1/2q \leq x \leq 1/2q+\beta/2) \\
f_4(x) & (1/2q+\beta/2 \leq x \leq 1/2q+\beta) \\
0 & (\text{otherwise}), 
\end{cases} 
\\
g_\beta(x)&=\chi_{[1/2q, 1/2q+\beta]}(x)\sqrt{f(x)-f(x)^2}.
\end{align*}
We note that, in either case, $f$ and $g$ are infinitely differentiable functions. Putting $e_\beta$ by 
\[ e_\beta=v^{-q'}g(u)+f(u)+g(u)v^{q'}, \]
where $f(u)$ and $g(u)$ belong to the Fr\'echet *-algebra $F^*(u)$ generated by $u$, we have the following lemma:
\begin{lemma}
$e_\beta$ cited above is a projection in $\tori$.
\end{lemma}

\begin{proof}
This follows from Connes \cite{co2}.
\end{proof}

 Another projection $e_{\beta'}$ is constructed by the similar way as $v$ and $u^{-1}$ in place of $u$ and $v$, and as $q'$ and $\beta'$ in place of $q$ and $\beta$ respectively. 

\begin{lemma}\label{th:2}
The projections $e_\beta, \alpha_{e^{\nipai p/q}, 1}(e_\beta), \dots , \alpha^{q-1}_{e^{\nipai p/q}, 1}(e_\beta)$ are mutually \\
orthogonal. So are the projections $e_{\beta'}, \alpha_{1, e^{-\nipai p'/q'}}(e_{\beta'}), \dots , \alpha^{q'-1}_{1, e^{-\nipai p'/q'}}(e_{\beta'})$.
\end{lemma}

\proof
We have that
\[ \alpha_{e^{\nipai p/q}, 1}(e_\beta)=
v^{-q'}g(e^{\nipai p/q}u)+f(e^{\nipai p/q}u)+g(e^{\nipai p/q}u)v^{q'}.
\]
Since the supports of $g$ and $g(e^{\nipai p/q}\cdot)$ are disjoint, we see for example that 
\begin{align*}
e_\beta\alpha_{e^{\nipai p/q}, 1}(e_\beta)&=
v^{-q'}g(u)v^{-q'}g(e^{\nipai p/q}u)+f(u)v^{-q'}g(e^{\nipai p/q}u) \\
&\quad +g(u)v^{q'}f(e^{\nipai p/q}u)+g(u)v^{q'}g(e^{\nipai p/q}u)v^{q'} \\
&=v^{-2q'}g(e^{-\nipai q'\theta}u)g(e^{\nipai p/q}u)+v^{q'}g(e^{\nipai q'\theta}u)f(e^{\nipai p/q}u) \\
&\quad +v^{-q'}f(e^{-\nipai p/q}u)+v^{q'}g(e^{\nipai q'\theta}u)g(e^{\nipai p/q}u)v^{q'} \\
&=v^{-2q'}g(e^{\nipai\beta}u)g(e^{\nipai p/q}u)+v^{-q'}f(e^{\nipai\beta}u)g(e^{\nipai p/q}u) \\
&\quad +v^{-q'}g(e^{-\nipai\beta}u)f(e^{\nipai p/q}u)+v^{q'}g(e^{-\nipai\beta}u)g(e^{\nipai p/q}u)v^{q'}. 
\end{align*}

When $\beta\geq 1/2q$, since $\supp f=[0, 1/q]$ and $\supp g=[\beta, 1/q]$, we have 
\begin{align*}
\supp g(e^{\nipai\beta}\cdot)&=[2\beta, 1/q+\beta], & \supp g(e^{-\nipai\beta}\cdot)&=[0, 1/q-\beta] \\
\supp g(e^{-\nipai p/q}\cdot)&=[\beta+p/q, (p+1)/q], & \supp f(e^{\nipai\beta}\cdot)&=[\beta, \beta+1/q] \\
\supp f(e^{\nipai p/q}\cdot)&=[p/q, (p+1)/q]. & 
\end{align*}
Using the fact that $p$ and $q$ are mutually prime, we conclude that the supports of $g(e^{\nipai \beta}\cdot)$ and $g(e^{\nipai p/q}\cdot)$ are disjoint and so on, which implies that $e_\beta\alpha_{e^{\nipai p/q}, 1}(e_\beta)=0$. By the analogous argument, we also have that the above equation holds when $\beta<1/2q$. By the same way, we see that
\[ \alpha^k_{e^{\nipai p/q}, 1}(e_\beta)\alpha^l_{e^{\nipai p/q}, 1}(e_\beta)=0 \] for $k, l \in \{0, 1, \cdots , q-1\}$ with $k\not=l$, as desired. Similarly, we can prove that the projections $e_{\beta'}, \alpha_{1, e^{-\nipai p'/q'}}(e_{\beta'}), \dots , \alpha^{q'-1}_{1, e^{-\nipai p'/q'}}(e_{\beta'})$ are also mutually orthogonal.
\qed

Now we define the elements $e_1$ and $e_2$ by
\[ e_1=\sum_{k=0}^{q'-1}(\alpha')^k(e_{\beta'}), \quad 
e_2=1-\sum_{k=0}^{q-1}\alpha^k(e_\beta), \]
where $\alpha=\alpha_{e^{\nipai p/q}, 1}, \alpha'=\alpha_{1, e^{-\nipai p'/q'}}$. By the previous proposition, both $e_1$ and $e_2$ are projections in $\tori$. Furthermore, we have that $\tau(e_\beta)=\beta, \tau(e_{\beta'})=\beta'$, where $\tau(x)$ is the canonical trace of $x\in T^2_\theta$.

\begin{lemma}\label{th:3}
The projections $e_1$ and $e_2$ are unitarily equivalent 
in $\tori$.
\end{lemma}

\proof
First of all, we show that $\tori$ is algebraically simple. Let $\ide$ be a non-zero *-ideal of $\tori$. Since the closure $\overline{\ide}$ of $\ide$ in $T^2_\theta$ is a closed *-ideal of $T^2_\theta$, it follows by the algebraic simplicity of $T^2_\theta$ that $\overline{\ide}$ must be equal to $T^2_\theta$. Then, there is an element $x\in\ide$ such that $\| 1-x\|<1$, so that the spectrum of $x$ does not include the origin of $\C$. Since the function $h(t)=1/t$ is holomorphic on the spectrum of $x$, it follows that $h(x)=x^{-1}\in \tori$.  Hence, $1=x^{-1}x\in\ide$, which implies that $\ide=T^2_\theta$, as claimed. 

Next, we have to verify that stable rank of $\tori$ is equal to one, i.e., the set of all invertible elements of $\tori$ is dense in $\tori$. If we would have this fact, $\tori$ has cancellation property (cf. Rieffel \cite{r2,r3}). Take any element $a\in \tori$. We may assume that $a\geq 0$. Then, for $\forall \e >0$, there exists an invertible element $b\geq 0$ in $T^2_\theta$ such that $\| a-b \|<\e/2$ (note that $T^2_\theta$ is of stable rank one.). By the density of $\tori$, we can find an element $c\in\tori$ with $c\geq 0$ and $\| b-c \| <\e/2$. We act $\tori$ on $L^2(T)$ defined before. Let us show that $c$ is invertible as an operator on $L^2(T)$. If $\xi\in \ker c$ and $\|b-c\|<\e/2$, we have 
\[ \|(b-c)\xi\|=\|b\xi\|<\frac{\e}{\,2\,}\|\xi\|. \]
Since $\e$ is arbitrary, we see that $\xi=0$, which means that $c$ is an injective operator. We note that we can find a positive number $\e/2>\delta>0$ such that $\|b\xi\|\geq \delta\|\xi\|$ for any $\xi \in L^2(T)$. We then have for any $\xi\in L^2(T)$,
\[ \|c\xi\|\geq |\|(b-c)\xi\|-\|b\xi\||\geq \left|\delta-\frac{\e}{\, 2\,}\right|\|\xi\|, \]
which implies that $c^{-1}$ is bounded. By triangle inequality, $\|a-c\|\leq \|a-b\|+\|b-c\|<\e$. Consequently, the stable rank of $\tori$ is one.

Now recall that $\tau(e_1)=\tau(e_2)$, we thus have $[e_1]=[e_2]\in K_0(\tori)$. Since $\tori$ has cancellation property, they are unitarily equivalent in $\tori$.
\qed

Let $\theta=[a_0, a_1, \dots , a_n, \dots ]$ be the continued fraction expansion and define the matrices $P_1, P_2, \cdots $ by
\[ P_n=
\begin{pmatrix}
a_{4n} & 1 \\
1 & 0
\end{pmatrix}
\begin{pmatrix}
a_{4n-1} & 1 \\
1 & 0
\end{pmatrix}\begin{pmatrix}
a_{4n-2} & 1 \\
1 & 0
\end{pmatrix}\begin{pmatrix}
a_{4n-3} & 1 \\
1 & 0
\end{pmatrix}
\]
for $n\geq 1$. Moreover, we put
\[ \begin{pmatrix}
q_{2n} \\
q_{2n-1}
\end{pmatrix}
=P_nP_{n-1}\cdots P_1
\begin{pmatrix}
1 \\ 0
\end{pmatrix}
\]
and 
\[ \fre=M_{q_{2n}}(C^\infty(T))\oplus M_{q_{2n-1}}(C^\infty(T)). \]
For each $n\geq 1$, we construct homomorphisms $\pi_n^\infty : \fre \to \mathfrak{A}_{n+1}$ as follows: we write $P_{n+1}=\left( \begin{smallmatrix} a&b \\ c&d \end{smallmatrix}\right)$. Let $z\in \ct$ be the canonical unitary generator of $\ct$. The element
\[ \begin{pmatrix}
z & & \\
& \ddots & \\
& & z 
\end{pmatrix}
\oplus O_{q_{2n-1}} \in \fre=M_{q_{2n}}(C^\infty(T))\oplus M_{q_{2n-1}}(C^\infty(T))
\]
should be mapped to the element
\[ \begin{pmatrix}
J_a & & & & & \\
& \ddots & & & & \\
& & J_a & & & \\
& & & O_b & & \\
& & & & \ddots & \\
& & & & & O_b
\end{pmatrix}
\oplus
\begin{pmatrix}
J_c' & & & & & \\
& \ddots & & & & \\
& & J_c' & & & \\
& & & O_d & & \\
& & & & \ddots & \\
& & & & & O_d
\end{pmatrix}
\in \mathfrak{A}_{n+1}
\]
$(=(\underbrace{J_a\oplus \cdots \oplus J_a}_{q_{2n}}\oplus \underbrace{O_b \oplus \cdots \oplus O_b}_{q_{2n-1}})\oplus (\underbrace{J_c'\oplus \cdots \oplus J_c'}_{q_{2n}}\oplus \underbrace{O_d \oplus \cdots \oplus O_d}_{q_{2n-1}}))$, where 
\[ J_k=\begin{pmatrix}
0 & & & z \\
1 & \ddots & & \\
& \ddots & \ddots & \\
& & 1 & 0 
\end{pmatrix},
\qquad
J_k'=\begin{pmatrix}
0 & & & 1 \\
1 & \ddots & & \\
& \ddots & \ddots & \\
& & 1 & 0 
\end{pmatrix}
\in M_k(\ct)
\]
and $O_l$ means the $l \times l$ zero matrix. Any element $(a_{ij})\oplus O_{q_{2n-1}} \in M_{q_{2n}}(\C)\oplus M_{q_{2n-1}}(\C)\subset \fre$ should be mapped to
\[ \begin{pmatrix}
a_{11}I_a & \cdots & a_{1q_{2n}}I_a & \\
\vdots & & \vdots & \\
a_{q_{2n}, 1}I_a & \cdots & a_{q_{2n}, q_{2n}}I_a  \\
& & & O_{bq_{2n-1}}
\end{pmatrix}
\oplus 
\begin{pmatrix}
a_{11}I_c & \cdots & a_{1q_{2n}}I_c & \\
\vdots & & \vdots & \\
a_{q_{2n}, 1}I_c & \cdots & a_{q_{2n}, q_{2n}}I_c  \\
& & & O_{dq_{2n-1}}
\end{pmatrix}
, \]
where $I_a, I_c$ are the $a\times a, c\times c$ identity matrices respectively. The second direct summand of $\fre$ should be mapped into $\mathfrak{A}_{n+1}$ by the similar way as $q_{2n}$ replaced by $q_{2n-1}$, $a$ and $c$ by $b$ and $d$ respectively, and interchanging the places to whose elements are mapped from upper left-hand side to lower right-hand side. It is easily verified that these $\pi_n^\infty$ are smooth inclusions.

Next, we need the following proposition. We define
\begin{align*}
 e_{kk}&=\alpha^{k-1}(e_\beta) & (k&=1, 2, \dots , q-1) \\
 \intertext{and}
 e_{kk}'&=(\alpha')^{k-1}(e_{\beta'}) & (k&=1, 2, \dots , q'-1).
\end{align*}

\begin{lemma}
Let $e_{22}ve_{11}=e_{21}|e_{22}ve_{11}|$ be the polar decomposition of $e_{22}ve_{11}$. Then, $e_{21}=e_{22}ve_{11}$.
\end{lemma}

\proof
We write $x=ve_{11}$. Since $x^*x=e_{11}v^*ve_{11}=e_{11}$, we have $|x|=e_{11}$. Thus, $x=ve_{11}$ is the polar decomposition of $x$, which implies that it is a surjective operator since $v$ is unitary. Hence, it follows that $\overline{\ran e_{22}}=\overline{\ran e_{22}ve_{11}}$, where $\overline{V}$ is the closure of a linear subspace $V$ of the Hilbert space $L^2(T)$. Furthermore, it is also verified that $\overline{\ran e_{11}}=\overline{\ran |e_{22}ve_{11}|}$. Note that $e_{22}ve_{11}=(e_{22}ve_{11})e_{11}$. By uniqueness of polar decomposition, we deduce that $e_{21}=e_{22}ve_{11}$, as desired.
\qed

By the similar way, we put $e_{21}'=e_{22}'ue_{11}'$. Our goal in this section is to construct the $F^*$-subalgebras generated by some unitaries, which is isomorphic to $M_{q_{2n}}(C^\infty(T))\oplus M_{q_{2n-1}}(C^\infty(T))$. For this, since $q_{2n-1}$ and $q_{2n}$ are mutually prime, we can find an integer $p_{2n-1}, p_{2n}$ with $\left(\begin{smallmatrix}
p_{2n-1} & p_{2n} \\ q_{2n-1} & q_{2n} \end{smallmatrix}\right)\in SL(2, \Z)$ and $p_n/q_n \to \theta$ as $n\to \infty$. With the same notations as above, we set
\[ \begin{pmatrix}
p' & p \\
q' & q
\end{pmatrix}
=
\begin{pmatrix}
p_{2n} & p_{2n-1} \\ q_{2n} & q_{2n-1} \end{pmatrix}
\]
and $\beta=\beta_n=p_{2n-1}-q_{2n-1}\theta, \beta'=\beta'_n=q_{2n}\theta-p_{2n}$, and so on. First of all, we check the following fact although it seems to be known:

\begin{lemma}\label{th:deri}
For arbitrary $h\in \ct$, $\delta_j(h(u))=h'(u)\delta_j(u) \, (j=1, 2)$, where $h'$ is the first derivative of $h$.
\end{lemma}

\proof
If $h(x)=\sum_{\nu=-m}^{n}a_\nu x^\nu$ is a Laurent polynomial, we have 
\begin{align*}
\delta_1(h(u))&=\delta_1\left(\sum_{\nu=-m}^{n}a_\nu u^\nu\right) =\sum_{\nu=-m}^{n}a_\nu \nu iu^\nu \\
&=\left(\sum_{\nu=-m}^{n}a_\nu \nu u^{\nu-1}\right)iu=h'(u)\delta_1(u).
\end{align*}
For any $h\in \ct$, we can find a family of Laurent polynomials $\{p_n\}_{n\geq 1}$ such that $p_n \to h$ with respect to the seminorms $\{\|\cdot \|_{k, l}\}$. For $m, n\geq 1$, we have
\begin{align*}
\delta_1(p_n(u)-p_m(u))&=(p_n'(u)-p_m'(u))\delta_1(u) \\
&=(p_n'(u)-p_m'(u))u.
\end{align*}
Since $\{p_n(u)\}_n$ is Cauchy, $\{\delta_1(p_n(u))\}_{n\geq 1}$ is also a Cauchy sequence. Using the fact that $\delta_1$ is a closed operator, we get
\begin{align*}
\delta_1(h(u))&=\lim_{n\to\infty}\delta_1(p_n(u)) \\
&=\lim_{n\to\infty}p_n'(u)\delta_1(u)=h'(u)\delta_1(u).
\end{align*}
As $\delta_2(u)=0$, it is clear that $\delta_2(h(u))=0=h'(u)\delta_2(u)$. This completes the proof.
\qed

In what follows, we use the notations $e_{11}^{(n)}=e_{\beta_n}, (e_{11}')^{(n)}=e'_{\beta_n}$ and so on for $n\geq 1$. Denoting $r_m=p_m/q_m$ for any integer $m\geq 1$, we define $u_n=u_{n, 1}+u_{n, 2}$ and $v_n=v_{n, 1}+v_{n, 2}$, where
\begin{align*}
u_{n, 1}&=\sum_{j=0}^{q_{2n}-1}e^{\nipai r_{2n}j}
\alpha^j_{e^{\nipai r_{2n}}, 1}(e_{11}^{(n)}), 
&\\
u_{n, 2}&=\sum_{j=0}^{q_{2n-1}-1}\alpha^j_{1, e^{-\nipai r_{2n-1}}}((e_{21}')^{(n)}) \\
v_{n, 1}&=\sum_{j=0}^{q_{2n}-1}\alpha^j_{e^{\nipai r_{2n}}, 1}(e_{21}^{(n)}), & \\
v_{n, 2}&=\sum_{j=0}^{q_{2n-1}-1}e^{-\nipai r_{2n-1}j}\alpha^j_{1, e^{-\nipai r_{2n-1}}}((e_{11}')^{(n)}).
\end{align*}
We note that since 
\begin{align*}
\alpha^{q_{2n}-1}(e_{21}^{(n)})&\in e_{11}^{(n)}\tori e_{q_{2n}q_{2n}}^{(n)} \\
(\alpha')^{q_{2n-1}-1}((e'_{21})^{(n)})&\in (e'_{11})^{(n)}\tori (e'_{q_{2n-1}q_{2n-1}})^{(n)},
\end{align*}
where $e_{q_{2n}q_{2n}}^{(n)}=\alpha^{q_{2n}-1}_{e^{\nipai r_{2n}}, 1}(e_{11}^{(n)})$ and $(e'_{q_{2n-1}q_{2n-1}})^{(n)}=\alpha^{q_{2n-1}-1}_{1, e^{-\nipai r_{2n-1}}}((e_{11}')^{(n)})$, we can find a unitary $v_{1q_{2n}}\in e_{11}^{(n)}\tori e_{11}^{(n)}$ (resp. $u'_{1q_{2n-1}}\in (e'_{11})^{(n)}\tori (e'_{11})^{(n)}$) such that $\alpha^{q_{2n}-1}(e_{21}^{(n)})=v_{1q_{2n}}e_{1q_{2n}}^{(n)}$ (resp. $(\alpha')^{q_{2n-1}-1}((e'_{21})^{(n)})=u'_{1q_{2n-1}}(e'_{1q_{2n-1}})^{(n)}$). By Lemma \ref{th:2}, we have 
\begin{align*}
u_{n, 1}u_{n, 1}^*
&=\left(\sum_{j=0}^{q_{2n}-1}e^{\nipai r_{2n}j}
\alpha^j_{e^{\nipai r_{2n}}, 1}(e_{11}^{(n)})\right) \\
&\quad \cdot \left(\sum_{j=0}^{q_{2n}-1}e^{-\nipai r_{2n}j}\alpha^j_{e^{\nipai r_{2n}}, 1}(e_{11}^{(n)})\right) \\
&=\sum_{j, m}e^{\nipai r_{2n}(j-m)}\alpha^j_{e^{\nipai r_{2n}}, 1}(e_{11}^{(n)})\alpha^m_{e^{\nipai r_{2n}}, 1}(e_{11}^{(n)}) \\
&=\sum_{j=0}^{q_{2n}-1}\alpha^j_{e^{\nipai r_{2n}}, 1}(e_{11}^{(n)})=1-e_2^{(n)}.
\end{align*}

Similarly, $u_{n, 1}^*u_{n, 1}=1-e_2^{(n)}, v_{n, 2}v_{n, 2}^*=v_{n, 2}^*v_{n, 2}=e_1^{(n)}$. Moreover, we have 
\begin{align*}
u_{n, 2}u_{n, 2}^*&=
\left(\sum_{j=0}^{q_{2n-1}-2}(e'_{2+j, 1+j})^{(n)}+u'_{1q_{2n-1}}(e'_{1q_{2n-1}})^{(n)}\right) \\
&\quad \cdot\left(\sum_{j=0}^{q_{2n-1}-2}(e'_{1+j, 2+j})^{(n)}+(e'_{q_{2n-1}1})^{(n)}(u'_{1q_{2n-1}})^*\right) \\
&=\left((e'_{21})^{(n)}+\dots +(e'_{q_{2n-1}, q_{2n-1}-1})^{(n)}\right) \\
&\quad \cdot \left((e'_{12})^{(n)}+\dots +(e'_{q_{2n-1}-1, q_{2n-1}})^{(n)}\right) \\
&\quad +\left((e'_{21})^{(n)}+\dots +(e'_{q_{2n-1}, q_{2n-1}-1})^{(n)}\right)u'_{1q_{2n-1}}(e'_{1q_{2n-1}})^{(n)} \\
&\quad +(e'_{q_{2n-1}, 1})^{(n)}u'_{1q_{2n-1}}\left((e'_{12})^{(n)}+\dots +(e'_{q_{2n-1}-1, q_{2n-1}})^{(n)}\right) \\
&\quad +(e'_{q_{2n-1}, 1})^{(n)}(u'_{1q_{2n-1}})^*u'_{1q_{2n-1}}(e'_{1q_{2n-1}})^{(n)},
\end{align*}
where 
\[ (e'_{k, k-1})^{(n)}=\alpha^{k-2}_{1, e^{-\nipai r_{2n-1}}}((e'_{11})^{(n)}), \quad (e_{k-1, k})^{(n)}=((e_{k, k-1})^{(n)})^* \]
 for $k=2, \dots , q_{2n-1}$. Since $u'_{1q_{2n-1}}$ is a unitary in $(e'_{11})^{(n)}\tori (e'_{11})^{(n)}$, it follows that the second and the third terms above are 0 and
\begin{align*}
(e'_{q_{2n-1}, 1})^{(n)}(u'_{1q_{2n-1}})^*u'_{1q_{2n-1}}(e'_{1q_{2n-1}})^{(n)}&=(e'_{q_{2n-1}1})^{(n)}(e'_{11})^{(n)}(e'_{1q_{2n-1}})^{(n)} \\
&=(e'_{q_{2n-1}q_{2n-1}})^{(n)}. 
\end{align*}
Thus we have
\[ u_{n, 2}u_{n, 2}^*=(e'_{11})^{(n)}+\dots +(e'_{q_{2n-1}-1, q_{2n-1}-1})^{(n)}+(e'_{q_{2n-1}q_{2n-1}})^{(n)}=e_1^{(n)}. \]
The same calculations show that
\[ u_{n, 2}^*u_{n, 2}=e_1^{(n)}, \quad v_{n, 1}v_{n, 1}^*=v_{n, 1}^*v_{n, 1}=1-e_2^{(n)}. \]
Moreover, we have
\begin{align*}
v_{n, 1}u_{n, 1} 
&=\left(e_{21}^{(n)}+\dots +e_{q_{2n}, q_{2n}-1}^{(n)}+u_{1q_{2n}}e_{1q_{2n}}^{(n)}\right)
\left(e_{11}^{(n)}+\dots +\omega^{q_{2n}-1}e_{q_{2n}q_{2n}}^{(n)}\right) \\
&=e_{21}^{(n)}+\dots +\omega^{q_{2n}-2}e_{q_{2n}q_{2n}-1}^{(n)}+\omega^{q_{2n}-1}u_{1q_{2n}}e_{1q_{2n}}^{(n)} \\
\intertext{and}
u_{n, 1}v_{n, 1}
&=\left(e_{11}^{(n)}+\dots +\omega^{q_{2n}-1}e_{q_{2n}q_{2n}}^{(n)}\right)\left(e_{21}^{(n)}+\dots +e_{q_{2n}, q_{2n}-1}^{(n)}+u_{1q_{2n}}e_{1q_{2n}}^{(n)}\right) \\
&=e_{11}^{(n)}u_{1q_{2n}}e_{1q_{2n}}^{(n)}+\omega e_{21}^{(n)}+\dots +\omega^{q_{2n}-1}e_{q_{2n}q_{2n}-1}^{(n)},
\end{align*}
where 
\begin{align*}
e_{kk}^{(n)}&=\alpha^{k-1}_{e^{\nipai r_{2n}}, 1}(e_{\beta_n}) \quad (k=2, \dots , q_{2n}-1), \\
e_{k, k-1}^{(n)}&=\alpha^{k-2}_{e^{\nipai r_{2n}}, 1}(e_{21}^{(n)}), \quad e_{k-1, k}^{(n)}=(e_{k, k-1}^{(n)})^*  \quad (k=2, \dots , q_{2n}) 
\end{align*}
and $\omega =e^{\nipai r_{2n}}$. Using the fact that $u_{1q_{2n}}\in e_{11}^{(n)}\tori e_{11}^{(n)}$ and $\omega^{q_{2n}}=1$, we have
\[ v_{n, 1}u_{n, 1}=e^{-\nipai r_{2n}}u_{n, 1}v_{n, 1}. \]
To sum up, we get the following:
\begin{lemma} 
The following hold:
\renewcommand{\labelenumi}{{\rm (\theenumi)}}

\begin{enumerate}
\item $u_{n,1}$ and $u_{n, 2}$ are unitaries in $(1-e_2^{(n)})\tori(1-e_2^{(n)})$ and so are $u_{n, 2}$ and $v_{n, 2}$ in $e_1^{(n)}\tori e_1^{(n)}$. 
\item $u_{n, 1}v_{n, 1}=e^{\nipai r_{2n}}v_{n, 1}u_{n, 1}, \quad u_{n, 2}v_{n, 2}=e^{\nipai r_{2n-1}}v_{n, 2}u_{n, 2}$.
\end{enumerate}
\end{lemma}

Now we construct subalgebras isomorphic to $M_{q_{2n}}(C^\infty (T))\oplus M_{q_{2n-1}}(C^\infty (T))$. Let $\{e_{ij}^{(n)}\}_{1\leq i, j\leq q_{2n}}$ be the matrix units constructed by 
\[ \{e_{11}^{(n)}, e_{22}^{(n)}, \dots e_{q_{2n}q_{2n}}^{(n)}, e_{21}^{(n)}, \dots , e_{q_{2n}, q_{2n}-1}^{(n)} \}. \]
 We then see the following lemma:

\begin{lemma}
The $F^*$-algebras $F^*(\{e_{ij}^{(n)}\}_{1\leq i, j\leq q_{2n}}, v_{1q_{2n}})$ generated by $\{e_{ij}^{(n)}\}_{1\leq i, j\leq q_{2n}}$ and $v_{1q_{2n}}$ are isomorphic to $M_{q_{2n}}(C^\infty (T))$ for all integers $n\geq 1$.
\end{lemma}

\proof
Consider the continuous field $S\ni t\mapsto e_{\beta_n}$ defined by Elliott and Evans \cite{ee}, where $S$ is a closed subinterval in $(0, \infty )$. The functions $f$ and $g$ appeared in the construction of $e_{\beta_{n}}$ are depend on $t\in S$, so that we write $f=f_t, \, g=g_t$. It is not difficult to verify that
\[ \|f_t^{(\nu )}-f_{t_0}^{(\nu )}\|_\infty \, , \,  \|g_t^{(\nu )}-g_{t_0}^{(\nu )}\|_\infty \to 0 \]
as $t\to t_0$ for any integer $\nu \geq 0$, where $f^{(\nu)}$ stands for the $\nu$-th derivatives of $f\in\ct$ and $\|\cdot \|_\infty$ is the supremum norm on $\ct$. Then our statement of this lemma follows immediately.
\qed

By the same way, it follows that the $F^*$-algebra $F^* (\{(e'_{ij})^{(n)}\}, u'_{1q_{2n-1}})$ generated by $\{(e'_{ij})^{(n)}\}_{1\leq i, j\leq q_{2n-1}}$ and $u'_{1q_{2n-1}}$ is isomorphic to $M_{q_{2n-1}}(\ct )$, where $\{(e'_{ij})^{(n)}\}_{1\leq i, j\leq q_{2n-1}}$ are the matrix units generated by
\[ \{ (e'_{11})^{(n)}, \dots (e'_{q_{2n-1}q_{2n-1}})^{(n)}, (e'_{21})^{(n)}, \dots (e'_{q_{2n-1}, \, q_{2n-1}-1})^{(n)}\}. \]

%%%%%%%%%%%%%%%%%%%%%%%%%%%%%%%%%%%%%%%%%%%%%%%%%%%%%%%%%%%%%%%
%\begin{align*}
%u_{n, 1}&\mapsto  
%\begin{pmatrix}
%z & & &  \\
% & \omega_{2n}z & & \\
% &  & \ddots & \\
% & & & \omega_{2n}^{q_{2n}-1}z
%\end{pmatrix}\oplus O , \\
%u_{n, 2}&\mapsto
%O\oplus  
%\begin{pmatrix}
%0 & & & z \\
%1 & \ddots & & \\
% & \ddots & \ddots & \\
% & & 1 & 0
%\end{pmatrix}, \\
%v_{n, 1}&\mapsto  
%\begin{pmatrix}
%0 & & & z \\
%1 & \ddots & & \\
% & \ddots & \ddots & \\
% & & 1 & 0
%\end{pmatrix}
%\oplus O
%\intertext{and}
%v_{n, 2}&\mapsto O\oplus 
%\begin{pmatrix}
%z & & &  \\
% & \omega_{2n-1}^{-1}z & & \\
% &  & \ddots & \\
% & & & \omega_{2n-1}^{-(q_{2n-1}-1)}z
%\end{pmatrix},
%\end{align*}
%where $\omega_n=e^{2\pi ir_n}$.
%%%%%%%%%%%%%%%%%%%%%%%%%%%%%%%%%%%%%%%%%%%%%%%%%%%%%%

\begin{lemma}\label{deri}
For each $h\in \ct$ and any integer $k\geq 1$, there exist $\{a_{\nu, k}\}\subset \R$ such that 
\[ \delta_1^k(h(u))=\sum_{\nu=1}^k a_{\nu, k}h^{(\nu)} (u)u^\nu
\quad  \quad (\nu=1, \dots , k). \]
\end{lemma}

\proof
For $k=1$, by Proposition \ref{th:deri}. If this statement holds for some $k\geq 1$, one has
\begin{align*}
\delta_1^{k+1}(h(u))&=\delta_1\left(\sum_{\nu=1}^k a_{\nu, k}h^{(\nu)}(u)u^\nu\right) \\
&=\sum_{\nu=1}^k a_{\nu, k}\delta_1(h^{(\nu)}(u)u^\nu) \\
&=\sum_{\nu=1}^k a_{\nu, k}\left( h^{(\nu+1)}(u)u\cdot u^\nu+i\nu h^{(\nu)}(u)u^\nu\right) \\
&=\sum_{\nu=1}^k a_{\nu, k}\left( h^{(\nu+1)}(u)u^{\nu+1}+i\nu h^{(\nu)}(u)u^\nu\right) \\
&=\sum_{\nu=2}^{k+1}a_{\nu -1, k}h^{(\nu)}(u)u^\nu
+\sum_{\nu=1}^k ia_{\nu, k}\nu h^{(\nu)}(u)u^\nu .
\end{align*}
Thus, we have 
\[ a_{\nu, k+1}=\sum_{\nu=2}^{k+1}a_{\nu -1, k}+
\sum_{\nu=1}^k ia_{\nu, k}\nu , \]
this ends the proof.

\qed

We note that the coefficients $a_{\nu, k}$ do not depend on the choice $h$.

By Lemma\,\ref{deri}, we have
\begin{align*}
\| \delta_1^k(f_n(u))-\delta_1^k(f_m(u)) \| 
&=\left\| \sum_{\nu=1}^k a_{\nu, k}\left(f_n^{(\nu)}(u)-f_m^{(\nu)}(u)\right)u^\nu\right\| \\
&\leq \sum_{\nu=1}^k |a_{\nu, k}|\|f_n^{(\nu)}(u)-f_m^{(\nu)}(u)\| \to 0 \quad
(n, m \to \infty),
\end{align*}
which means that $\{\delta_1^k(f_n(u))\}_n$ is a Cauchy sequence. Analogously, we see that $\{\delta_1^k(g_n(u))\}_n$ is also Cauchy.

By construction, the following fact follows:

\begin{lemma}
Let $F^* (u_n, v_n)$ be the $F^*$-algebras generated by $u_n$ and $v_n$. Then, they are equal to $F^*(\{e_{ij}^{(n)}\}, v_{1q_{2n}})\oplus F^*(\{(e'_{ij})^{(n)}\}, u'_{1q_{2n-1}})$.
\end{lemma}

\begin{proof}
Since $u_{n, j}$ and $v_{n, j}$ $(j=1, 2)$ are all periodic unitaries, their spectra are finite. Then the projections appeared in the spectral decompositions of $u_{n, j}, v_{n, j}$ are unitarily equivalent to $e_{ij}^{(n)}$s by the properties that $F^*(u_{n, j})$ and $F^*(v_{n, j})$ are closed under the holomorphic functional calculus.
\end{proof}

\begin{lemma}
For any integers $k, l\geq 0$, 
\[ \lim_{n\to \infty}\|u-u_n\|_{k, l}=\lim_{n\to\infty}\|v-v_n\|_{k, l}=0. \]
\end{lemma}

\proof
At first, we have to verify that the sequence $\{\delta_1^k(e_{\beta_n})\}_n$ is Cauchy. By construction of $e_{\beta_n}$, we have, for $n, m\geq 1$,
\begin{align*}
\|\delta_1^k(e_{\beta_n})-\delta_1^k(e_{\beta_m})\|
&\leq \|\delta_1^k(v^{-q_{2n-1}}g_n(u)-v^{-q_{2m-1}}g_m(u))\| \\
&\quad +\| \delta_1^k(f_n(u)-f_m(u)) \|+\|\delta_1^k(g_n(u)v^{q_{2n-1}}-g_m(u)v^{q_{2m-1}})\| \\
&=\|v^{-q_{2n-1}}\delta_1^k(g_n(u))-v^{-q_{2m-1}}\delta_1^k(g_m(u))\| \\
&\quad +\|\delta_1^k(f_n(u))-\delta_1^k(f_m(u))\| \\
&\quad +\|\delta_1^k(g_n(u))v^{q_{2n-1}}-\delta_1^k(g_m(u))v^{q_{2m-1}}\|.
\end{align*}

Since $p_{2n-1}/q_{2n-1} \to \theta$, the last term of the above calculation tends to $0$ as $n, m\to \infty$. Therefore, $\{\delta_1^k\circ\delta_2^l(u(1-e_2^{(n)})-u_{n, 1})\}_n$ is Cauchy. Similarily, the sequence $\{\delta_1^k\circ\delta_2^l(ue_1^{(n)}-u_{n, 2})\}_n$ is also a Cauchy sequence. Hence, by \cite{r}, 
\[ u(1-e_2^{(n)})-u_{n, 1} \to 0, \quad ue_1^{(n)}-u_{n, 2} \to 0 \]
as $n \to \infty$. Using the fact that $\delta_1^k\circ\delta_2^l$ are closed, the sequences above tend to $0$ as $n\to \infty$. Consequently, 
\begin{align*}
\|u-u_n\|_{k, l}&\leq \|u(1-e_2^{(n)})-u_{n, 1}\|_{k, l}+\|ue_1^{(n)}-u_{n, 2}\|_{k, l} \\
&\to 0 \qquad (n\to \infty).
\end{align*}
By the similar argument, we have $\|v-v_n\|_{k, l}\to 0$ as $n\to\infty$, this ends the proof.
\qed

Combining all together in this section, we conclude that our key fact follows:

\begin{proposition}
Given an irrational number $\theta \in (0, 1)$, $\tori$ is isomorphic to the Fr\'echet $^*$-inductive limit
\[ \varinjlim (M_{q_{2n}}(\ct)\oplus M_{q_{2n-1}}(\ct), \pi_n^\infty). \]
\end{proposition}

\section{Entire Cyclic Cohomology of Fr\'echet Inductive Limits}\label{fre}

Let $\{\fre, i_n \}_{n\geq 1}$ be a family of Fr\'echet *-algebras and $i_n : \fre \to \mathfrak{A}_{n+1}$ Fr\'echet *-imbeddings. We can form the Fr\'echet *-inductive limit $\varinjlim \mathfrak{A}_n$, which is denoted by $\mathfrak{A}$. In this section, we prove that the projective limit $\varprojlim H_\e^*(\fre)$ of the entire cyclic cohomologies $\varprojlim H_\e^*(\fre)$ is isomorphic to $H_\e^*(\ore)$. Let $[\,\cdot\,]_{\fre}$ be the entire cyclic cohomology classes on $\fre$, and the maps $\widehat{i_n}^* : H_\e^{{\rm ev}}(\ore_{n+1})\to H_\e^{{\rm ev}}(\fre)$ are defined by
\[ \widehat{i_n}^*([(\varphi_{2k}^{(n+1)})_k]_{\ore_{n+1}})=[(i_n^{\otimes (2k+1)})^*\varphi_{2k}^{(n+1)}]_{\fre},
\]
where
\[ (i_n^{\otimes (2k+1)})^*\varphi_{2k}^{(n+1)}(a_0, \dots , a_{2k})
=\varphi_{2k}^{(n+1)}(i_n(a_0), \dots , i_n(a_{2k})) \]
for $\retsu \in \fre$.
First of all, we define the notion of projective limit as follows:

\begin{definition} \label{teigi}
The projective limit $\entirea$ of $H_\e^{{\rm ev}}(\fre)$ is the space of sequences $\{ [(\varphi_{2k}^{(n)})_k ]_{\fre} \}_n \in \prod_{n\geq 1} H_{\varepsilon}^{\rm{ev}}(\fre)$ such that for any $n\geq 1$,
\[ \widehat{i_n}^*([(\varphi^{(n+1)}_{2k})_k]_{\ore_{n+1}})= [(\varphi^{(n)}_{2k})_k]_{\fre} \]
with the property that for any $k\geq 0, l\geq 1$,
\[ \sup_{n\geq 1}\| \varphi^{(n)}_{2k} \|_l < \infty, \]
where
\[ \| \varphi_{2k}^{(n)} \|_l = \sup_{a_j \in \fre, \, \| a_j \|_l\leq 1} | \varphi_{2k}^{(n)}(a_0, \dots , a_{2k}) |. \]
We define $\varprojlim H_\e^{{\rm od}}(\fre)$ in the similar way as in the even case. $\{ [(\varphi_{2k}^{(n)})_k ]_{\fre} \}_n = \{ [(\psi_{2k}^{(n)})_k ]_{\fre} \}_n$ if and only if there exists $\{ [(\theta_{2k+1}^{(n)})_k ]_{\fre} \}_n \in \varprojlim H_\e^{{\rm od}}(\fre)$ such that
\[ \varphi_{2k}^{(n)}-\psi_{2k}^{(n)}=b\theta_{2k-1}^{(n)}+B\theta_{2k+1}^{(n)} \]
for any $n\geq 1, k\geq 0$.
\end{definition}

Let us construct two maps between $\entirea$ and $H_\e^{{\rm ev}}(\ore)$. First of all, we define $\Phi : \entireb \to \entirea$ by
\[ \Phi([(\varphi_{2k})_k]_\ore)=\{[(\varphi_{2k}|_{\fre})_k]_{\fre}\}_n, \]
where $[\,\cdot\,]_\ore$ means the same symbol as $[\,\cdot\,]_{\fre}$. Actually it is well-defined. In fact, if $[(\varphi_{2k})_k]_\ore=[(\varphi'_{2k})_k]_\ore$ then there exists an odd entire cyclic cocycle $\theta=(\theta_{2k+1})_k$ such that $(\varphi_{2k}-\varphi'_{2k})_k=(b+B)(\theta_{2k+1})_k,$ where $b+B$ is the derivation on entire cyclic cocycles. It is trivial that $(\varphi_{2k}|_{\fre}-\varphi'_{2k}|_{\fre})_k=(b+B)(\theta_{2k+1}|_{\fre})_k$ for each integer $n\geq 1$. This means that $\{[(\varphi_{2k}|_{\fre})_k]_{\fre}\}_n=\{[(\varphi'_{2k}|_{\fre})_k]_{\fre}\}_n$. Moreover, 
\[ \sup_{n\geq 1} \| \varphi_{2k}^{(n)}|_{\fre} \|_l=\| \varphi_{2k} \|_l < \infty, \]
which implies $[(\varphi_{2k}^{(n)}|_{\fre})_k]_{\fre}\in H_\e^{{\rm ev}}(\fre)$. \\

 Now we construct the inverse map $\Psi$ of $\Phi$. For any $\{[(\varphi_{2k}^{(n)})_k]_{\fre}\}_n$ $\in \entirea$ and $\retsu \in \ore$, we can take sequences $\{b_j^{(m)}\}_m$ for $j=0, \dots , 2k$ which converge to $a_j$ as $m\to\infty$ with respect to the seminorms $\|\cdot\|_l$ on $\varinjlim \fre$. Choose integers $N(m)\geq 1$ such that $b_j^{(m)}\in \mathfrak{A}_{N(m)}$ for any $0\leq j\leq 2k$. We may assume that $N(m)=m$ by taking a larger number between $N(m)$ and $m$. We have that for $m>m'$, there exists an odd entire cocycle $\theta^{(m')}=(\theta^{(m')}_{2k+1})_k$ on $\mathfrak{A}_{m'}$ such that
\begin{align}\label{shiki1}
&\varphi_{2k}^{(m)}(b_0^{(m)}, \dots , b_{2k}^{(m)})-\varphi_{2k}^{(m')}(b_0^{(m')}, \dots , b_{2k}^{(m')}) \\
&\quad =(b\theta_{2k-1}^{(m')}+B\theta_{2k+1}^{(m')})(b_0^{(m')}, \dots , b_{2k}^{(m')}). \notag
\end{align}
By Hahn-Banach theorem, we can extend $\varphi_{2k}^{(m)}$ and $\varphi_{2k}^{(m')}$ to $\widetilde{\varphi_{2k}}^{(m)}$ and $\widetilde{\varphi_{2k}}^{(m')}$ on $\ore$ such that
\[ \|\widetilde{\varphi_{2k}}^{(m)}\|_l=\|\varphi_{2k}^{(m)}\|_l, \quad \|\widetilde{\varphi_{2k}}^{(m')}\|_l=\|\varphi_{2k}^{(m')}\|_l \]
for any $l\geq 1$.
\begin{lemma}
For any $a_0, \dots , a_{2k} \in \ore$, the sequence
\[ \{ \widetilde{\varphi_{2k}}^{(m)} (\retsu)\}_m \]
is bounded.
\end{lemma}

\proof
We have
\begin{align*}
|\widetilde{\varphi_{2k}}^{(m)}(\retsu)|&\leq |\widetilde{\varphi_{2k}}^{(m)}(a_0-b_0^{(m)}, a_1, \dots , a_{2k})| \\
&+|\widetilde{\varphi_{2k}}^{(m)}(b_0^{(m)}, a_1-b_1^{(m)}, a_2, \dots , a_{2k})| \\
&+\cdots \\
&+|\widetilde{\varphi_{2k}}^{(m)}(b_0^{(m)}, \dots , b_{2k-1}^{(m)}, a_{2k}-b_{2k}^{(m)})| \\
&+|\widetilde{\varphi_{2k}}^{(m)}(b_0^{(m)}, \dots , b_{2k}^{(m)})|.
\end{align*}
By the above equation (\ref{shiki1}), 
\begin{align*}
&\widetilde{\varphi_{2k}}^{(m)}(b_0^{(m)}, \dots , b_{2k}^{(m)}) \\
&=\varphi_{2k}^{(m)}(b_0^{(m)}, \dots , b_{2k}^{(m)}) \\
&=\varphi_{2k}^{(m')}(b_0^{(m')}, \dots , b_{2k}^{(m')})
+(b\theta_{2k-1}^{(m')}+B\theta_{2k+1}^{(m')})(b_0^{(m')}, \dots , b_{2k}^{(m')})
\end{align*}
is a constant independent of $m$. Using the hypothesis in Definition \ref{teigi} and Hahn-Banach theorem, it follows that $\lim_{m\to \infty}|\widetilde{\varphi_{2k}}^{(m)}(\retsu)|$ is dominated by the constant $|\varphi_{2k}^{(m')}(b_0^{(m')}, \dots , b_{2k}^{(m')})+(b\theta_{2k-1}^{(m')}+B\theta_{2k+1}^{(m')})(b_0^{(m')}, \dots , b_{2k}^{(m')})|$. In particular, the sequence $\{|\widetilde{\varphi_{2k}}^{(m)}(\retsu)|\}_m$ is bounded.
\qed \\

Therefore, by taking the subsequence of $\{|\widetilde{\varphi_{2k}}^{(N)}(\retsu)|\}_N$, we may assume that
\[ \lim_{N\to \infty}\widetilde{\varphi_{2k}}^{(N)}(\retsu) \]
exists, so that we define
\[ \widetilde{\varphi_{2k}}(\retsu)=\lim_{N\to \infty}\widetilde{\varphi_{2k}}^{(N)}(\retsu). \]
Here we note that
\[ \widetilde{\varphi_{2k}}(\retsu)=\lim_{m\to \infty}\widetilde{\varphi_{2k}}^{(m)}(b_0^{(m)}, \dots , b_{2k}^{(m)}). \]
In fact, by the same reason as before, we have
\begin{align*}
&|\widetilde{\varphi_{2k}}^{(m)}(\retsu)-\widetilde{\varphi_{2k}}^{(m)}(b_0^{(m)}, \dots , b_{2k}^{(m)})| \\
&\leq |\widetilde{\varphi_{2k}}^{(m)}(a_0-b_0^{(m)}, a_1, \dots , a_{2k})| \\
&\quad +\cdots \\
&\quad +|\widetilde{\varphi_{2k}}^{(m)}(b_0^{(m)}, \dots , b_{2k-1}^{(m)}, a_{2k}-b_{2k}^{(m)})| \to 0
\end{align*}
as $m \to \infty$. Using the above preparation, we shall show the following fact:

\begin{lemma}
$(\widetilde{\varphi_{2k}})_k$ is an entire cyclic cocycle on $\ore$.
\end{lemma}

\proof
Let $\Sigma$ be a bounded subset of $\ore$ and $\retsu \in \Sigma$. Then we can choose sequences $\{b_j^{(m)}\}_m \subset \bigcup \fre$ for $j=0, \dots , 2k$ such that $b_j^{(m)}\to a_j$ as $m\to \infty$ with respect to the topology induced by the seminorms $\| \cdot \|_l$ on $\ore$. In this case, the set
\[ \Sigma_0 = \{b_j^{(m)}\in \bigcup \fre \, | \, j=0, \dots , 2k, m\in \mathbb{N} \} \]
is bounded in $\ore$. So, by the equation (\ref{shiki1}),
\begin{align*}
|\widetilde{\varphi_{2k}}(\retsu)|&=\lim_{m\to\infty}|\widetilde{\varphi_{2k}}^{(m)}(b_0^{(m)}, \dots , b_{2k}^{(m)})| \\
&\leq |\widetilde{\varphi_{2k}}^{(1)}(b_0^{(1)}, \dots , b_{2k}^{(1)})| \\
&\quad +|(b\theta_{2k-1}^{(1)}+B\theta_{2k+1}^{(1)})(b_0^{(1)}, \dots , b_{2k}^{(1)})|.
\end{align*}
As $(\varphi_{2k}^{(1)})_k$ and $(b\theta_{2k-1}^{(1)}+B\theta_{2k+1}^{(1)})_k$ are entire on $\mathfrak{A}_1$, 
\[ |\widetilde{\varphi_{2k}}(\retsu)|\leq Ck! \]
for some constant $C>0$ independent of $m$, which implies that $(\widetilde{\varphi_{2k}})_k$ is entire.
\qed \\

Now we are ready to define a map $\Psi : \varprojlim H_\e^{{\rm ev}}(\fre) \to H_\e^{{\rm ev}}(\ore)$ in the following fashion:
\[ \Psi(\{ [(\varphi_{2k}^{(n)})_k ]_{\fre} \}_n)=[(\widetilde{\varphi_{2k}})_k]_\ore . \]
We have to verify that the definition is well-defined. Let
\[ \{ [(\varphi_{2k}^{(n)})_k ]_{\fre} \}_n = \{ [(\psi_{2k}^{(n)})_k ]_{\fre} \}_n \in \varprojlim H_\e^{{\rm ev}}(\fre) . \]
Then for any $n\geq 1$, there exists an odd entire cyclic cocycles $\theta^{(n)}=(\theta_{2k+1}^{(n)})_k$ on $\fre$ such that
\[ \varphi_{2k}^{(n)}(b_0, \dots , b_{2k})-\psi_{2k}^{(n)}(b_0, \dots , b_{2k})
=(b\theta_{2k-1}^{(n)}+B\theta_{2k+1}^{(n)})(b_0, \dots b_{2k}) \]
for $b_0, \dots , b_{2k}\in \fre.$ By the above argument, there exists an odd entrie cyclic cocycle $\widetilde{\theta}=(\widetilde{\theta_{2k+1}})_k$ on $\ore$. Then by the definition of $b+B$, we have that
\begin{align*}
&(b\theta_{2k-1}^{(n)}+B\theta_{2k+1}^{(n)})(\retsu) \\
&=\lim_{m\to\infty}(b\theta_{2k-1}^{(m)}+B\theta_{2k+1}^{(m)})(b_0^{(m)}, \dots , b_{2k}^{(m)}) \\
&=\lim_{m\to\infty}\left(\widetilde{\varphi_{2k}}^{(m)}(b_0^{(m)}, \dots , b_{2k}^{(m)})-\widetilde{\psi_{2k}}^{(m)}(b_0^{(m)}, \dots , b_{2k}^{(m)})\right) \\
&=\widetilde{\varphi_{2k}}(\retsu)-\widetilde{\psi_{2k}}(\retsu),
\end{align*}
which implies that $[(\widetilde{\varphi_{2k}})_k]_\ore=[(\widetilde{\psi_{2k}})_k]_\ore$.
\begin{proposition}\label{th:10}
The following isomorphism holds as a vector space over $\C$:
\[ \varprojlim H_\e^*(\fre) \simeq H_\e^{*}(\ore). \]
\end{proposition}

\proof
We prove just in the even case.
For any $[(\varphi_{2k})_k]_\ore \in H_\e^{{\rm ev}}(\ore)$, we have
\[ \Psi \circ \Phi ([(\varphi_{2k})_k]_\ore)=\Psi(\{[(\varphi_{2k}|_{\fre})_k]_{\fre}\}_n)=[(\widetilde{\varphi_{2k}|_{\fre}})_k]_\ore . \]

\noindent
For any $\retsu \in \ore$, we take sequences $\{b_j^{(m)}\}_{m}$ \, $(j=0, \cdots, 2k)$ which converge to $a_j$ as $m\to\infty$ and $b_j^{(m)}\in \mathfrak{A}_{m}$ for $j=0, \cdots , 2k$. Then,
\begin{align*}
\widetilde{\varphi_{2k}|_{\fre}}(\retsu)
&= \lim_{m \to \infty} \varphi_{2k}|_{\mathfrak{A}_{m}}(b_0^{(m)}, \dots , b_{2k}^{(m)}) \\
&=\varphi_{2k}(\retsu).
\end{align*}
This implies that $\widetilde{\varphi_{2k}|_{\mathfrak{A}_n}}=\varphi_{2k}$ , which means that $\Psi \circ \Phi$ is the identity on $\entireb$. On the other hand, for any $\{[(\varphi_{2k}^{(n)})_k]_{\fre}\}_n \in \entirea$, we have
\[ \Phi\circ\Psi(\{[(\varphi_{2k}^{(n)})_k]_{\fre}\}_n)=\Phi([(\widetilde{\varphi_{2k}})_k]_\ore )=\{[(\widetilde{\varphi_{2k}}|_{\fre})_k]_{\fre}\}_n. \]
Since for $b_0, \dots , b_{2k} \in \fre$, we have
\begin{align*}
\widetilde{\varphi_{2k}}|_{\fre}(b_0, \dots , b_{2k})&=\lim_{m\to\infty}\widetilde{\varphi_{2k}}^{(m)}(b_0, \dots , b_{2k}) \\
&=\lim_{m\to\infty}\varphi_{2k}^{(m)}(b_0, \dots , b_{2k}) \\
&=\varphi_{2k}^{(n)}(b_0, \dots , b_{2k}).
\end{align*}
Thus $\Phi\circ\Psi(\{[(\varphi_{2k}^{(n)})_k]_{\fre}\}_n)=\{[(\varphi_{2k}^{(n)})_k]_{\fre}\}_n$.
 Hence $\Phi\circ\Psi$ is also the identity on $\entirea$. Therefore, the proof is completed. 
\qed \\

\noindent
{\it{Remark}}. \, We here prefer the original defintion by Connes \cite{co} to prove our main result although Meyer \cite{meyer} obtained the above Proposition by means of analytic cyclic theory.

\section{Entire Cyclic Cohomology of $\tori$}\label{main}

Summing up the argument discussed in the previous sections, we are ready to obtain the next main result

\begin{theorem}
The entire cyclic cohomology $H_\e^*(\tori)$ of the noncommutative $2$-torus $\tori$ is isomorphic to $\C^4$ as linear spaces, especially
\[
\begin{cases}
H_\e^{\rm{ev}}(\tori)=HP^{\rm{ev}}(\tori)\simeq \C^2 \\
H_\e^{\rm{od}}(\tori)=HP^{\rm{od}}(\tori)\simeq \C^2,
\end{cases}
\]
where $HP^*(\tori)$ is the periodic cyclic cohomology of $\tori$.
\end{theorem}

\proof
By Lemma~\ref{th:10}, we have
\begin{align*}
H_\e^*(\tori)&\simeq H_\e^*(\varinjlim (\ct\otimes (M_{q_{2n}}(\C)\oplus M_{q_{2n-1}}(\C)), \
\pi_n^\infty)) \\
&\simeq \varprojlim H_\e^*((\ct\otimes (M_{q_{2n}}(\C)\oplus M_{q_{2n-1}}(\C)), (\pi^\infty_n)^* ) 
\end{align*}

We have the following decomposition by applying Khalkhali \cite{k}'s Proposition 7 in the case of $F^*$-algebras:
\begin{align*}
 H_\e^*(C^\infty (T)\otimes (M_{q_{2n}}(\C)\oplus M_{q_{2n-1}}(\C)))
\simeq H_\e^*(C^\infty (T)\otimes (M_{q_{2n}}(\C)\oplus H_\e^*(C^\infty (T)\otimes M_{q_{2n-1}}(\C))) 
\end{align*}
\noindent
We also deduce applying Khalkhali \cite{k}'s Theorem 6 in the case of $F^*$-algebras that 
\[ H_\e^*(C^\infty (T)\otimes (M_{q}(\C)) \simeq H_\e^*(C^\infty (T))~~(q \geq 1) \]
\noindent
Since the above two phenomena are shown for $HP^*(\tori)$ as well and we can see that 
\[H_\e^j(C^\infty (T))=HP^j(C^\infty (T)) \simeq \C ~~(j=\rm{ev,od})\]
\noindent
(Connes \cite{co},Thm 2(page 208) and Thm 25(page 382)), then we obtain that 
\[H_\e^j(C^\infty (T)\otimes (M_{q}(\C)) \simeq HP^j(C^\infty (T)\otimes (M_{q}(\C)) ~~~~(j=\rm{ev,od})\]
\noindent
We then have the the following commutative diagram :
\[
\begin{CD}
HP^{\rm{ev}}(\mathfrak{A}_{n+1})@>\simeq>i^*> H_\e^{\rm{ev}}(\mathfrak{A}_{n+1}) \\
@V(\pi^\infty_n)^*VV  @VV(\pi^\infty_n)^*V \\
HP^{\rm{ev}}(\mathfrak{A}_{n}) @>i^*>\simeq> H_\e^{\rm{ev}}(\mathfrak{A}_{n})
\end{CD}
\]
where $i^*$ is the canonical inclusion map.
Then we work on the periodic cyclic cohomology in what follows: we consider homomorphisms
\begin{align*}
(\pi^\infty_n)^* : HP^{\rm{ev}}(C^\infty (T)\otimes (M_{q_{2n+2}}(\C)\oplus M_{q_{2n+1}}(\C))) \\
\to
HP^{\rm{ev}}(C^\infty (T)\otimes (M_{q_{2n}}(\C)\oplus M_{q_{2n-1}}(\C))). 
\end{align*}
Now we note that 
\begin{align*}
&\, HP^{\rm{ev}}(C^\infty (T)\otimes (M_{q_{2n+2}}(\C)\oplus M_{q_{2n+1}}(\C))) \\
&\simeq HP^{\rm{ev}}(C^\infty (T)\otimes M_{q_{2n+2}}(\C))
\oplus HP^{\rm{ev}}(C^\infty (T)\otimes M_{q_{2n+1}}(\C))
\end{align*}
and moreover, we have seen that
\begin{align*}
HP^{\rm{ev}}(C^\infty (T)\otimes M_{q}(\C))
&\simeq HP^{\rm{ev}}(C^\infty (T))\otimes HP^{\rm{ev}}(M_q(\C)) \\
&\simeq \C \left[\int_T \right] \otimes \C \left[{\rm Tr}_q\right] \\
&\simeq \C \left[\int_T \otimes {\rm Tr}_q \right],
\end{align*}
where $\int_T$ and ${\rm Tr}_q$ are the usual integral on $C^\infty(T)$ and the trace on $M_q(\C)$ respectively. Here, we consider the following diagram:
\[
\begin{CD}
HP^{\rm{ev}}(\mathfrak{A}_{n+1})@>\simeq>> \C\left[\int_T\otimes {\rm Tr}_{q_{2n+2}}\right]\oplus \C\left[\int_T\otimes {\rm Tr}_{q_{2n+1}}\right] \\
@V(\pi^\infty_n)^*VV  @VV(\pi^\infty_n)^*V \\
HP^{\rm{ev}}(\mathfrak{A}_{n}) @>>\simeq> \C\left[\int_T\otimes {\rm Tr}_{q_{2n}}\right]\oplus \C\left[\int_T\otimes {\rm Tr}_{q_{2n-1}}\right],
\end{CD}
\]
where the horizonal isomorphisms are defined by
\begin{align*}
 HP^{\rm{ev}}(\mathfrak{A}_n) \to \C \left[\int_T\otimes {\rm Tr}_{q_{2n}}\right]\oplus\left[\int_T\otimes {\rm Tr}_{q_{2n-1}}\right] \\
\varphi\mapsto \varphi|_{(C^\infty(T)\otimes M_{q_{2n}}(\C))\oplus 0} \oplus \varphi|_{0\oplus (C^\infty(T)\otimes M_{q_{2n-1}}(\C))}. 
\end{align*}
We check that the diagram above is also commutative.

So, we regard $(\pi^\infty_n)^*$ as the linear map from $\C [\int_T \otimes {\rm Tr}_{q_{2n+2}}]\oplus \C [\int_T \otimes {\rm Tr}_{q_{2n+1}}]$ into $\C [\int_T \otimes {\rm Tr}_{q_{2n}}]\oplus \C [\int_T \otimes {\rm Tr}_{q_{2n-1}}]$.
Let us recall that we write the matrix $P_{n+1}$ by $\left( \begin{smallmatrix} a&b \\ c&d \end{smallmatrix}\right)$ used in the definition of $\pi^\infty_n$. Then we have
\begin{align*}
\left(\left(\int_T \otimes {\rm Tr}_{q_{2n+2}}\right)\oplus 0\right)(\pi^\infty_n (\xi))
&=a \left( \int_T \otimes {\rm Tr}_{q_{2n}}\right)({\bf 1}\otimes (x_{ij})) \tag{$*$}\label{con} \\
&\quad  +b \left( \int_T \otimes {\rm Tr}_{q_{2n-1}}\right)({\bf 1}\otimes (y_{ij}))  
\end{align*}
for each $\xi=({\bf 1}\otimes (x_{ij}))\oplus ({\bf 1}\otimes (y_{ij}))
\in (C^\infty(T)\otimes M_{q_{2n}}(\C))\oplus (C^\infty(T)\otimes M_{q_{2n-1}}(\C))$, where ${\bf 1}$ is the function which evaluates 1 at each point of $T$. In fact, by the definition of $\pi^\infty_n$, we have
\begin{align*}
\pi^\infty_n (({\bf 1}\otimes (x_{ij})\oplus ({\bf 1}\otimes (y_{ij}))
&=
\begin{pmatrix}
x_{11}I_a & \dots & x_{1q'}I_a & & & \\
\vdots & & \vdots & & & \\
x_{q'1}I_a & \dots & x_{q'q'}I_a & & & \\
& & & y_{11}I_b & \dots & y_{1q}I_b \\
& & & \vdots & & \vdots \\
& & & y_{q1}I_b & \dots & y_{qq}I_b
\end{pmatrix}
\\
&\oplus
\begin{pmatrix}
x_{11}I_c & \dots & x_{1q'}I_c & & & \\
\vdots & & \vdots & & & \\
x_{q'1}I_c & \dots & x_{q'q'}I_c & & & \\
& & & y_{11}I_d & \dots & y_{1q}I_d \\
& & & \vdots & & \vdots \\
& & & y_{q1}I_d & \dots & y_{qq}I_d
\end{pmatrix},
\end{align*}
where $q=q_{2n-1}, \, q'=q_{2n}$ and so on. Then, it follows that
\begin{align*}
\left(\left(\int_T \otimes {\rm Tr}_{q_{2n+2}}\right)\oplus 0 \right)(\pi^\infty_n(\xi)) 
&=a\sum_{i=1}^{q_{2n}}x_{ii} + b\sum_{i=1}^{q_{2n-1}}y_{ii} \\
&=a \left(\int_T \otimes {\rm Tr}_{q_{2n}}\right)({\bf 1}\otimes (x_{ij})) \\&\quad +b \left(\int_T \otimes {\rm Tr}_{q_{2n-1}}\right)({\bf 1}\otimes (y_{ij})).
\end{align*}
Similarly, we have 
\begin{align*}
\left(0\oplus \left(\int_T \otimes {\rm Tr}_{q_{2n+1}}\right)\right)(\pi^\infty_n (\xi))
&=c \left( \int_T \otimes {\rm Tr}_{q_{2n}}\right)({\bf 1}\otimes (x_{ij})) \tag{$**$}\label{con2} \\
&\quad  +d \left( \int_T \otimes {\rm Tr}_{q_{2n-1}}\right)({\bf 1}\otimes (y_{ij})).  
\end{align*}
On the other hand, we check that
\begin{align*}
\left(\left(\int_T\otimes {\rm Tr}_{q_{2n+2}}\right)\oplus 0\right)(\pi^\infty_n((z^k\otimes I_{q_{2n}})\oplus 0))&=0 \\
\left(0\oplus\left(\int_T\otimes {\rm Tr}_{q_{2n+1}}\right)\right)(\pi^\infty_n((z^k\otimes I_{q_{2n}})\oplus 0))&=0 \\
\left(\left(\int_T\otimes {\rm Tr}_{q_{2n+2}}\right)\oplus 0\right)(\pi^\infty_n(0\oplus (z^k\otimes I_{q_{2n-1}})))&=0 \\
\intertext{and} 
\left(0\oplus\left(\int_T\otimes {\rm Tr}_{q_{2n+1}}\right)\right)(\pi^\infty_n(0\oplus (z^k\otimes I_{q_{2n-1}})))&=0 
\end{align*}
for each integer $k\geq 1$. Indeed, for example, it is easily verified that if \[ 
\begin{pmatrix}
0 & & & z \\
1 & \ddots &  &  \\
& \ddots & \ddots & \\
& & 1 & 0
\end{pmatrix}
\in M_q (C^\infty (T)),
\]
\begin{align*}
\begin{pmatrix}
0 & & & z \\
1 & \ddots &  &  \\
& \ddots & \ddots & \\
& & 1 & 0
\end{pmatrix}^k =
\begin{cases}
z^\nu\otimes I_q & (k\equiv 0 \mod q) \\
\begin{pmatrix}
0 & &  {\large *} \\
& \ddots &  \\
{\large *} & &  0
\end{pmatrix}
& (k\not\equiv 0 \mod q)
\end{cases}
\end{align*}
for some integer $\nu\geq 1$. Thus, we have that
\begin{align*}
\left(\int_T\otimes {\rm{Tr}}_q\right)
\left(
\begin{pmatrix}
0 & & & z \\
1 & \ddots &  &  \\
& \ddots & \ddots & \\
& & 1 & 0
\end{pmatrix}^k
\right)&=
\begin{cases}
{\displaystyle \int_T} z^\nu dz & (k\equiv 0 \mod q) \\
0 & (k\not\equiv 0 \mod q)
\end{cases}
\\
&=0.
\end{align*}

Since the space of Laurent polynomials are dense in $C^\infty(T)$ with respect to Fr\'echet topology, we then conclude that (\ref{con}) and (\ref{con2}) hold for every $\xi\in \mathfrak{A}_n$. Hence, it is verified that $(\pi^\infty_n)^*$ is an isomorphism by the fact that
\begin{align*}
\det
\begin{pmatrix}
a & b \\
c & d
\end{pmatrix}
&=\det P_{n+1} \\
&=\det
\begin{pmatrix}
a_{4n+4} & 1 \\
1 & 0
\end{pmatrix}
\begin{pmatrix}
a_{4n+3} & 1 \\
1 & 0
\end{pmatrix}\begin{pmatrix}
a_{4n+2} & 1 \\
1 & 0
\end{pmatrix}\begin{pmatrix}
a_{4n+1} & 1 \\
1 & 0
\end{pmatrix}
\\
&=1\neq 0
\end{align*}

Finally, we conclude that
\begin{align*}
H_\e^{\rm{ev}}(\tori)&\simeq \varprojlim (\C\oplus\C, (\pi^\infty_n)^*) \\
&\simeq \C^2.
\end{align*}
Analogously, the same consequence is obtained in the odd case. We note that
\begin{align*}
&\, HP^{\rm{od}}(C^\infty(T)\otimes M_q(\C)) \\
&\simeq HP^{\rm{od}}(C^\infty(T))\otimes HP^{\rm{ev}}(M_q(\C)) \\
&\simeq \C\left[\psi\otimes {\rm Tr}_q\right],
\end{align*}
where $\psi(f, g)=\int_T f(t)g'(t)dt$ for $f, g\in C^\infty (T)$.
 This ends the proof.
\qed

\section*{acknowledgements}
I would like to thank my supervisor Professor H. Takai for suggesting this problem and many useful advices. I am also very grateful to Professor A. Connes for his valuable suggestions.

\end{document}